\documentclass[a4paper,english,12pt]{article}
\usepackage[logonly]{trace}
\usepackage{babel}
\usepackage{amsfonts, amsmath, amssymb}
\usepackage{graphicx}
\usepackage{pstricks}
\usepackage{psfig}
\usepackage{multirow}
\usepackage{fancyhdr}
\usepackage{vmargin, fancybox}

\newcommand{\mathsym}[1]{{}}

\def\N{\textrm{I\kern-0.21emN}}

\def\R{\textrm{I\kern-0.21emR}}
\def\Q{\textrm{l\kern-0.5emQ}}

\begin{document}

\title{A note on solitary waves solutions of classical wave equations}

\author{Claire David\footnotemark[2] \footnote{Corresponding author:
david@lmm.jussieu.fr; fax number: (+33) 1.44.27.52.59.}   \\ \small{Universit\'e Pierre et Marie Curie-Paris 6}  \\
\footnotemark[2] \small{ Institut Jean Le Rond d'Alembert, UMR CNRS 7190,} \\
\small{Bo\^ite courrier $n^0 162$, 4 place Jussieu, 75252 Paris,
cedex 05, France}}


\maketitle

\begin{abstract}
The goal of this work is to determine whole classes of solitary
wave solutions general for wave equations.
\end{abstract}

\section{Introduction}
\label{sec:intro}

\noindent Consider a differential equation of the form:

\begin{equation}
\label{EqInitiale}  {\mathcal {F}} (u,\frac{\partial u^{r
}}{\partial x^r },\frac{\partial u^{ s}}{\partial  t^s}) = 0,
\end{equation}

\noindent The determination of travelling wave solutions of
specific cases of (\ref{EqInitiale}), such as the Burgers or
Burgers-Korteweg-De Vries equations, for instance, has been a
major topic in the past few years, and play a crucial role in the
study of wave equations. We presently aim at extending previous
results, and give the whole classes of solitary wave solutions
general for (\ref{EqInitiale}).

\noindent The paper is organized as follows. The general method is
exposed in Section \ref{Solitons}. A specific case is studied in
section \ref{Solutions}.

\section{Solitary waves}
\label{Solitons}

\noindent Following Feng \cite{feng1} and our previous work
\cite{David}, in which travelling wave solutions of the
\textit{CBKDV} equation were exhibited as combinations of
bell-profile waves and kink-profile waves, we aim at determining
travelling wave solutions of (\ref{EqInitiale}) (see
\cite{li},\cite{whitham},
\cite{ablowitz}, \cite{dodd}, \cite{johnson}, \cite{ince}, \cite{zhang3}, \cite{birk}, \cite{Polyanin}).\\
\noindent Following \cite{feng1}, we assume that equation
(\ref{EqInitiale}) has  travelling wave solutions of the form
\begin{equation}\label{ChgtVar} {u}({x}, {t}) ={u}(\xi), \quad \xi= {x}-v\,{t} \end{equation}
where $v$ is the wave velocity. Substituting (\ref{ChgtVar}) into
equation (\ref{EqInitiale}) leads to:

\begin{equation}
\label{EqXi}   {\mathcal {F}} ( {u}, {u}^{(r)},(-v)^s\, {u}^{(s)})
= 0,
\end{equation}

\noindent Performing an integration of (\ref{EqXi}) with respect
to $\xi$ leads to an equation of the form:

\begin{equation}
\label{EqXiInt}   {\mathcal {F}}_\xi^{\mathcal P} ( {u},
{u}^{(r)},(-v)^s\, {u}^{(s)}) = C,
\end{equation}

\noindent where $C$ is an arbitrary integration constant, which
will be the starting point for the determination of solitary waves
solutions.\\

\noindent In the previous works, this integration constant is
usually taken equal to zero. Yet, it should not be so, since it
can lead to a loss of solutions, as we are going to show it in the
following.

\section{Travelling Solitary Waves}

\label{Solutions}

\subsection{Hyperbolic Ansatz}

\noindent The discussion in the preceding section provides us
useful information when we construct travelling solitary wave
solutions for equation (\ref{EqXi}). Based on these results, in
this section, a class of travelling wave solutions is searched as
a combination of bell-profile waves and kink-profile waves of the
form
\begin{equation}
\label{sol} \tilde{u}(\tilde{x}, \tilde{t}) = \sum_{i = 1}^n \left
(U_i\; \text{tanh}^i \left [\,C_i (\tilde{x}-v\,\tilde{t})\,
\right ] + V_i \; \text{sech}^i \left
[\,C_i(\tilde{x}-v\,\tilde{t}+x_0)\, \right ] \right )+V_0
\end{equation}
where the $U_i's$, $V_i's$,  $C_i's$, $(i=1,\ \cdots,\ n)$, $V_0$
and $v$ are constants to be determined.
\\
\noindent In the following, $c$ is taken equal to 1.

\subsection{Theoretical analysis}

\noindent Substitution of (\ref{sol}) into equation
(\ref{EqXiInt}) leads to an equation of the form

\begin{equation}
\label{EqgenInt}
 \sum_{i,\, j, \,k} A_i\,\text{tanh}^i \big (C_i\,
\xi \big )\,\text{sech}^j \big (C_i \,\xi \big )\,\text{sinh}^k
\big (C_i\, \xi  \big )  =C
\end{equation}

\noindent the $A_i$ being real constants.\\
 \noindent The difficulty
for solving equation (\ref{EqgenInt}) lies in finding the values
of the constants $U_i$, $V_i$,  $C_i$, $V_0$ and $v$ by solving
the over-determined algebraic equations. Following \cite{feng1},
after balancing the higher-order derivative term and the leading
nonlinear term, we deduce $n=1$.\\

\noindent Then, following \cite{David} we replace
$\mbox{sech}({C_1} \,\xi)$ by $ \frac{2}{ e^{\,{C_1} \,\xi}+e^{\,-
{C_1} \,\xi } }$, $\mbox{sinh(}{C_1} \,\xi)$ by $ \frac {
e^{\,{C_1} \,\xi}-e^{\,- {C_1} \,\xi } }{2}$, $\mbox{tanh}({C_1}
\,\xi)$ by $\frac{e^{\, {C_1} \,\xi } -
      e^{\,- {C_1} \,\xi}}{e^{ \,{C_1} \,\xi   } + e^{\,-{C_1} \,\xi }}$, and multiply both sides by
      $(1+ e^{ 2\,\xi \, {C_1}})^2$, so that equation
(\ref{EqgenInt}) can be rewritten in the following form:
\begin {equation}
\label{Syst} \sum_{k=0}^{4} P_k ( U_1,\ V_1,\ C_1,\ v,\, V_0 )
\,e^{\,k \,C_1\, \xi} \; = \;0,
\end{equation}
where the $P_k$ $(k=0,\,...,\,4)$, are polynomials of $U_1$,
$V_1$, $C_1$, $V_0$  and $v$.

\noindent Depending wether (\ref{EqgenInt}) admits or no
consistent solutions,  spurious solitary waves solutions may, or
not, appear.

\subsection{A specific case}

\noindent Consider the specific case when (\ref{EqInitiale}) is
the equivalent equation of a \emph{DRP} scheme, the coefficients
of which will be denoted by $\gamma_{k}$, $k\in \{-m,m\}$ (see
\cite{ Chaos2}):

\begin{equation}
\scriptsize  \label{EqEq} \begin{aligned} &  -  {u}_{ {t}}
-\frac{\sigma}{2  } \,  {u}_{{  {t}}{  {t}}}
  +\frac{ 2\,\sigma  }{\mu\,Re_h} \,\displaystyle \sum _{k=1}^m
k\,\gamma_{k}  \,  {u}_{{  {x}}  }   =0
\end{aligned}
\end{equation}

\noindent

\noindent where $ Re_h$ denotes the mesh Reynolds number,
$\sigma$, the $cfl$ coefficient, and $\mu$, the viscosity.

\noindent Equation (\ref{EqXi}) is then given by:

\begin{equation}
\label{EqXi1} \scriptsize  \begin{aligned} &  -v\,\tilde{u}'( \xi
) -\frac{v^2\,\sigma}{2 } \,\tilde {u}''(\xi)
  +\frac{ 2\,\sigma  }{\mu\,Re_h} \,\displaystyle \sum _{k=1}^m
k\,\gamma_{k}  \,\tilde {u}'(\xi) =0
\end{aligned}
\end{equation}
\noindent Performing an integration of (\ref{EqXi1}) with respect
to $\xi$ yields:

\begin{equation}
\label{EqXi2}  \scriptsize  \begin{aligned} & -v\,\tilde{u} ( \xi
) -\frac{v^2\,\sigma}{2 }\,\tilde {u}' (\xi)
  +\frac{ 2\,\sigma  }{\mu\,Re_h} \,\displaystyle \sum _{k=1}^m
k\,\gamma_{k}   \,\tilde {u} (\xi) =C
\end{aligned}
\end{equation}

\noindent i. e.:

\begin{equation}
\label{EqXi3}  \scriptsize    \left \lbrace \frac{ 2\,\sigma
}{\mu\,Re_h} \,\displaystyle \sum _{k=1}^m k\,\gamma_{k}  -v
\right \rbrace \,\tilde {u} (\xi) -\frac{v^2\,\sigma}{2 }\,\tilde
{u}' (\xi)=C
\end{equation}

\noindent where $C$ is an arbitrary integration constant.

\noindent Substitution of (\ref{sol}) for $n=1$ into equation
(\ref{EqXi3}) leads to:

\begin{equation}
\label{EqXi4}  \scriptsize    \left \lbrace \frac{ 2\,\sigma
}{\mu\,Re_h} \,\displaystyle \sum _{k=1}^m k\,\gamma_{k}  -v
\right \rbrace \,  \left \lbrace U_1\, \text{tanh}  \left [\,C_1
\,\xi\, \right ] + V_1 \, \text{sech} \left [\,C_1\,\xi\, \right ]
+V_0  \right \rbrace -\frac{v^2\,\sigma}{2 }\, \left \lbrace U_1\,
\text{sech}^2  \left [\,C_1\,\xi\, \right ] - V_1 \, \frac {
\text{sinh} \left [\,C_1\,\xi\, \right ] }{ \text{cosh}^2 \left
[\,C_1\,\xi\, \right ] } \right \rbrace=C
\end{equation}

\noindent i. e.:

\begin{equation}
\label{EqXi4}  \scriptsize    \left \lbrace \frac{ 2\,\sigma
}{\mu\,Re_h} \,\displaystyle \sum _{k=1}^mk\,\gamma_{k}  -v \right
\rbrace \,  \left \lbrace U_1\,
\frac{e^{C_1\,\xi}-e^{-C_1\,\xi}}{e^{C_1\,\xi}+e^{-C_1\,\xi}}\, +
  \frac{2\,V_1}{e^{C_1\,\xi}+e^{-C_1\,\xi}} +V_0  \right
\rbrace -\frac{v^2\,\sigma}{2 }\, \left \lbrace U_1\, \left (
\frac{2}{e^{C_1\,\xi}+e^{-C_1\,\xi}} \right )^2 - 2\,V_1 \,
\frac{e^{C_1\,\xi}-e^{-C_1\,\xi}}{\left (
{e^{\,C_1\,\xi}+e^{\,-\,C_1\,\xi}} \right )^2} \right \rbrace=C
\end{equation}

\noindent Multiplying both sides by $\left ( {1+e^{\,2\,C_1\,\xi}}
\right )^2$ yields:

\begin{equation}
\label{EqXi4}  \scriptsize
\begin{aligned} &  \left \lbrace \frac{ 2\,\sigma
}{\mu\,Re_h} \,\displaystyle \sum _{k=1}^mk\,\gamma_{k}  -v \right
\rbrace \,  \left \lbrace U_1\,
  \left ({e^{\,4\,C_1\,\xi}-1} \right ) \, +
   2\,V_1\,\left (e^{\,3\,C_1\,\xi}+e^{\,C_1\,\xi} \right ) +V_0 \,\left ( {1+e^{\,2\,C_1\,\xi}}
\right )^2 \right \rbrace\\
&-\frac{v^2\,C_1\,\sigma \,C_1}{2 }\, \left \lbrace 4\, U_1  -
2\,V_1 \,\left ({e^{\,3\,C_1\,\xi}-1}\right ) \right \rbrace=C
\end{aligned}
\end{equation}

\noindent which is a fourth-order equation in $e^{\,C_1\,\xi}$.
This equation being satisfied for any real value of $\xi$, one
therefore deduces that the coefficients of $e^{\,k\,C_1\,\xi}$,
$k=0,\, \ldots, \,4$ must be equal to zero, i.e.:

\begin{equation}
\label{EqXi4}  \scriptsize \left \lbrace
\begin{aligned}
 &   2\,\left \lbrace \frac{ 2\,\sigma
}{\mu\,Re_h} \,\displaystyle \sum _{k=1}^mk\,\gamma_{k}  -v \right
\rbrace \,  \left \lbrace -U_1 +V_0   \right
\rbrace-\frac{v^2\,C_1\,\sigma}{2 }\, \left \lbrace 4\, U_1  +
2\,V_1
  \right \rbrace=C\\
 &  \left
\lbrace \frac{ 2\,\sigma }{\mu\,Re_h} \,\displaystyle \sum
_{k=1}^mk\,\gamma_{k}  -v \right \rbrace \,
   2\,V_1     =0\\
&   2\,\left \lbrace \frac{ 2\,\sigma }{\mu\,Re_h} \,\displaystyle
\sum _{k=1}^mk\,\gamma_{k}  -v \right \rbrace \,
 V_0=0\\
&  2\, \left \lbrace \frac{ 2\,\sigma }{\mu\,Re_h} \,\displaystyle
\sum _{k=1}^m k\,\gamma_{k}  -v \right \rbrace \,
   V_1
 + v^2\,C_1\,\sigma \,V_1   =0\\
&  \left \lbrace \frac{ 2\,\sigma }{\mu\,Re_h} \,\displaystyle
\sum _{k=1}^mk\,\gamma_{k} -v \right \rbrace \, \left \lbrace U_1
     +V_0 \,
 \right \rbrace=0\\
\end{aligned} \right.
\end{equation}

\noindent $v=\frac{ 2\,\sigma }{\mu\,Re_h} \,\displaystyle \sum
_{k=1}^mk\,\gamma_{k} \, , \, V_1 \neq 0$ leads to the trivial
null solution. Therefore, $V_1$ is necessarily equal to zero,
which implies:

\begin{equation}
\label{EqXi4}  \scriptsize \left \lbrace
\begin{aligned}
&v=\frac{ 2\,\sigma }{\mu\,Re_h} \,\displaystyle \sum
_{k=1}^mk\,\gamma_{k} \\
 &     U_1
   =-\frac{C}{2\,C_1\, v^2\,\sigma}\\
&  V_0 \in \R  \,\,\, , \,\,\,C_1 \in \R
\end{aligned} \right.
\end{equation}

\noindent It is easy to note that, if the integration constant $C$
had been taken equal to zero, the solitary waves of the considered
equation would have been loss.

\normalsize

\section{Conclusions}

The importance of choosing an integration constant which is not
equal to zero, in the determination of solitary wave solutions of
wave equations, has been carried out. We show that taking this
constant equal to zero leads to a loss of solutions.

\addcontentsline{toc}{section}{\numberline{}References}

\end{document}